\documentclass{amsart}

\usepackage{amsmath, enumerate, amsfonts, amssymb, amsthm, graphics, graphicx, hyperref}

\theoremstyle{definition}
\newtheorem{thm}{Theorem}[section]
\newtheorem{prop}[thm]{Proposition}
\newtheorem{cor}[thm]{Corollary}
\newtheorem{lem}[thm]{Lemma}
\theoremstyle{definition}

\theoremstyle{definition}
\numberwithin{equation}{section}
\newtheorem{example}[thm]{Example}

\def\ie{{i.e.,}\ }

\newcommand{\nor}{\trianglelefteq}
\newcommand{\sn}{\triangleleft \triangleleft \ }

\newcommand{\iso}{\cong}
\def\ord#1{\vert #1 \vert}
\def\ind#1#2{\vert #1\, :\,#2\vert}

\def\Z#1{\textrm{Z}(#1)}

\def\CD#1{\mathcal{C}\mathcal{D} (#1)}
\def\CL#1{\mathcal{C}\mathcal{L} (#1)}
\def\M#1{\mathfrak{M}_{#1}}

\title{\bf Some groups with computable Chermak-Delgado lattices}

\author{Ben Brewster}
\address{Ben Brewster, Department of Mathematical Sciences, Binghamton University \\
		Binghamton, New York 13902}
\email{ben@math.binghamton.edu}

\author{Elizabeth Wilcox}
\address{Elizabeth Wilcox, Department of Mathematics, Colgate University \\
   Hamilton, New York 13346}
\email{ewilcox@colgate.edu}

\date{\today}

\begin{document}


\markboth{\textsc{Computable Chermak-Delgado lattices}}{\textsc{Computable Chermak-Delgado lattices}}

\begin{abstract}
Let $G$ be a finite group and let $H \leq G$.  We refer to $\ord {H} \ord {C_G(H)}$ as the {\it Chermak-Delgado measure of $H$} with respect to $G$.  Originally described by A. Chermak and A. Delgado, the collection of all subgroups of $G$ with maximal Chermak-Delgado measure, denoted $\CD {G}$, is a sublattice of the lattice of all subgroups of $G$.  In this paper we note that if $H \in \CD G$ then $H$ is subnormal in $G$ and prove if $K$ is a second finite group then $\CD {G \times K} = \CD G \times \CD K$.  We additionally describe the $\CD {G \wr C_p}$ where $G$ has a non-trivial center and $p$ is an odd prime and determine conditions for a wreath product to be a member of its own Chermak-Delgado lattice.  We also examine the behavior of centrally large subgroups, a subset of the Chermak-Delgado lattice.
\end{abstract}

\maketitle

A. Chermak and A. Delgado \cite{cd} defined a family of functions from the set of subgroups of a finite group into the set of positive integers. Chermak and Delgado then used these functions to obtain a variety of results, including a proof that every finite group $G$ has a characteristic abelian subgroup $N$ such that $\ind {G}{N} \leq {\ind {G}{A}}^2$ for all abelian $A \leq G$.  

In \cite{Isaacs}, I. Martin Isaacs focused on one member of this family, which he referred to as {\it the Chermak-Delgado measure}.  Isaacs showed for a fixed group $G$ that the subgroups with maximal measure form a sublattice within the lattice of subgroups of $G$, which he referred to as the {\it Chermak-Delgado lattice} of $G$.  After observing a paucity of groups which were members of their own Chermak-Delgado lattice, it seemed natural to investigate their existence. Thus, in this paper we study the  Chermak-Delgado lattice of direct products and wreath products. We prove that its members are always subnormal in $G$ and find special conditions in which $G \wr H$ is in its own Chermak-Delgado lattice.  As a by-product of our efforts, we show that every 2-group can be embedded as a subnormal subgroup of a group that is a member of its Chermak-Delgado lattice.

Moreover, in a recent article G. Glauberman studies some large subgroups of the Chermak-Delgado lattice \cite{glauberman2006}.  We show that this collection of subgroups behaves nicely in direct products and wreath products $G \wr C_p$, where $C_p$ is the cyclic group of odd order $p$.

Throughout the paper we use the following familiar notation. For $n$ a positive integer we use $S_n$ to denote the symmetric group on $n$ points and $A_n$ to denote the alternating subgroup of $S_n$.  We use $C_n$, $D_n$, and $Q_n$ to represent the cyclic, dihedral, and quaternion group of order $n$ (respectively, and for applicable values of $n$).  If $D$ is a direct product with $G$ as one of its factors then $\pi_G$ will represent the natural projection map from $D$ onto  $G$.  If $D$ is the direct product of multiple copies of $G$ with itself, then we use $G_i$ to represent the $i^{\textrm{th}}$ factor in $D$ and $\pi_i$ to represent the projection map from $D$ onto $G_i$.

\section{Preliminaries}

Define the {\it Chermak-Delgado measure} of a subgroup $H$ with respect to a finite group $G$ with $H \leq G$ as
$$
m_G (H) = \ord H \ord {C_G (H)}.
$$
From the definition, it's clear that the groups discussed in this paper are necessarily finite.  The next two lemmas are straightforward to prove using just the definition of $m_G(H)$ and recollections about centralizers from introductory group theory courses.

\begin{lem}\label{lattice1} If $H \leq G$ then $m_G(H) \leq m_G (C_G(H))$, and if the measures are equal then $H = C_G(C_G(H))$. \end{lem}

\begin{lem}\label{lattice2} If $H, K \leq G$ then $m_G(H)m_G(K) \leq m_G(\langle H,K \rangle) m_G (H \cap K)$.  Moreover, equality occurs if and only if $\langle H, K \rangle = HK$ and $C_G (H \cap K) = C_G(H)C_G(K)$. \end{lem}

The full details of the proofs of these lemmas can be found in \cite[Section 1.G]{Isaacs}.  For any finite group $G$, let $\M G$ denote the maximal measure over all subgroups in $G$ and let the set of all subgroups $H \leq G$ with $m_G(H) = \M G$ be denoted by $\CD {G}$.  From Lemmas~\ref{lattice1} and \ref{lattice2} we see:

\begin{thm}\label{lattice} For a finite group $G$ the set $\CD {G}$ is a sublattice within the lattice of subgroups of $G$ and for all $H$, $K$ in $\CD {G}$ we have $\langle H, K \rangle = HK$.  Moreover, if $H \in \CD {G}$ then $C_G (H) \in \CD {G}$ and $H = C_G (C_G (H))$.  \end{thm}

The lattice described in Theorem~\ref{lattice} will be referred to as {\it the Chermak-Delgado lattice of $G$}.  Clearly $\CD {G}$ is a sublattice within in the lattice of subgroups of $G$.  For large or complex groups $G$, it can be a challenge to determine $\CD {G}$ by hand. The calculations for small groups and abelian groups is, on the other hand, refreshingly easy:

\begin{enumerate}[1. ]
\item Let $G$ be abelian. If $H \leq G$ then $m_G(H) = \ord H \ord G$.  Therefore the only subgroup of maximal measure is $G$ and $\CD {G} = \{ G \}$.  
\item Let $G = S_4$; then $m_G(G) = 24 = m_G (\Z {G})$.  With a little work, one can show that the measure of any other subgroup of $S_4$ is less than 24 -- for example $m_G(A_4) = 16$.  Hence $\CD {S_4} = \{ S_4, 1 \}$.
\item Let $G = S_3$; then $m_G(G) = m_G (\Z {G}) = 6$.  On the other hand, the subgroup $A_3$ is abelian and is also its own centralizer.  Thus $m_G(A_3) = 9$.  The subgroups of order 2 in $G$ are also their own centralizers, therefore these subgroups have measure 4.  Hence $\CD {G} = \{A_3\}$.
\item Consider $D_8$, the dihedral group of order 8.  There are 5 subgroups of $D_8$ with measure 16, these are: $D_8$, $\Z {D_8}$, and the three subgroups of order 4.  All other subgroups have smaller measure.  Hence $\CD {D_8}$ is also the lattice of normal subgroups of $D_8$.  
\item One can also show that $\CD{Q_8}$ is isomorphic to the lattice of normal subgroups of $Q_8$, which happens to be isomorphic to the lattice of normal subgroups of $D_8$.
\end{enumerate}

Observe, since both $D_8$ and $S_3$ can be represented as subgroups of $S_4$, that there is not a straightforward relation between $\CD U$ and $\CD G$ when $U \leq G$.  Of course, one notices that if $U \leq G$ then
$$
\begin{array}{rl}
\M U &= m_U (V) \qquad \textrm{ for some } V \leq U\\
&= \ord V \ord {C_U(V)}\\
&\leq \ord V \ord {C_G (V)}\\
&\leq m_G (V)\\
&\leq \M G.\\
\end{array}
$$

The next result is due to Wielandt, and can be found as Theorem 2.9 in \cite{Isaacs}.  Isaacs refers to the result as a ``Zipper Lemma'' and we continue that reference here.

\begin{thm}[``Zipper Lemma''] Suppose that $S \leq G$ where $G$ is a finite group and assume that $S \sn H$ for every proper subgroup $H$ of $G$ that contains $S$.  If $S$ is not subnormal in $G$ then there is a unique maximal subgroup of $G$ that contains $S$. \end{thm}

The Zipper Lemma makes way for the use of induction with regards to the Chermak-Delgado lattice. Another important fact regarding $U \in \CD G$ is the following:

\begin{prop}\label{useful}
Let $U \in \CD {G}$ for a finite group $G$.  If $S < G$ with both $U \leq S$ and $U C_G (U) \leq S$ then $U \in \CD S$.  \end{prop}

Proposition~\ref{useful} is easy to see -- when $U C_G(U) \leq S$ then $C_G (U) = C_S (U)$.  In fact, not only is $U \in \CD S$ but also $C_G(U)$ and $U C_G(U)$ are in $\CD S$.  This useful proposition, together with the Zipper Lemma, is enough to prove the following result.  

\begin{thm}\label{subnormal} Let $G$ be a group.  If $U \in \CD G$ then $U \sn G$. \end{thm}

\begin{proof} Assume that for every proper subgroup $U$ of $G$ if $X \in \CD U$ then $X \sn U$.  Now let $U \in \CD G$.  We show $V = U C_G(U) \sn G$, which is sufficient for the theorem since $U \nor V$.

If $V = G$ our conclusion holds, so assume that $V < G$.  For every $S < G$ with $V \leq S$ we know that $U \leq S$ and therefore $V \in \CD S$ by Proposition~\ref{useful}.  By induction $V \sn S$.  If $V$ is not subnormal in $G$ then we may apply the Zipper Lemma, resulting in the existence of a unique maximal subgroup $M$ of $G$ that contains $V$.  Notice, by the previous few sentences, that $V \sn M$.

Let $x \in G$; since $V^x \in \CD G$ we have $VV^x \in \CD G$ as well.  If $VV^x = G$ then there exists $v, v_0 \in V$ such that $vv_0^x = x$, and careful multiplication shows that $x = v_0v \in V$.  Thus if $VV^x = G$ then $V = G$.  We assumed $V < G$, though, thus $VV^x < G$.  There exists a proper maximal subgroup $N$ of $G$ that contains $VV^x$, however since $V < N$ and $M$ is the unique maximal subgroup containing $V$ we know $N = M$.

One can repeat the use of the Zipper Lemma on $V^x$ and determine that $M^x$ is the unique maximal subgroup of $G$ containing $V^x$.  Yet $M$ contains $V^x$, thus $M = M^x$ for all $x \in G$.  Hence $M \nor G$.  Subnormality is transitive, therefore $V \sn G$ as desired. \end{proof}

A trivial consequence of Theorem~\ref{subnormal} is that the Chermak-Delgado lattice of any simple group $S$ is $\{\Z S, S\}$ (of course $\Z S = S$ when $S$ is abelian).  Another easy consequence of Theorem~\ref{subnormal} is the expansion of Examples 2 and 3: Given a symmetric group $S_n$ for $n \geq 5$ we know that the only possible subgroups in $\CD {S_n}$ are $1$, $A_n$, and $S_n$.  Since the measure of $A_n$ will be less than that of $S_n$, we know that $\CD {S_n} = \{1, S_n\}$; therefore the Chermak-Delgado lattice of any symmetric group is completely determined.  
  
One might question whether Theorem~\ref{subnormal} can be strengthened, \ie are all subgroups in $\CD G$ actually normal in $G$?  The answer, demonstrated by the next example, is negative.

\begin{example}  Let $G$ be as follows.
$$
\begin{array}{rl}
G = \langle a, b, c, d \quad \mid & a^4 = b^2 = c^2 = d^2 = [a,b] = [b,c] = [b,d]\\
&  = [c, d] = [a, c]b = [a, d]c = 1 \rangle
\end{array}
$$
This presentation is convenient for computations, though $G$ actually is a 2-generator group.  A few calculations (made by hand or with GAP \cite{gap}) show that $X = \langle a, b \rangle$ is a member of $\CD G$.  One can show that $d$ does not normalize $X$, therefore $X \sn G$ with defect greater than 1.  There are a few other subgroups in $\CD {G}$ that are not normal, such as $\langle b, da \rangle$ and $\langle b, da^3 \rangle$, though showing by hand that these subgroups are not normal is tedious. \end{example}

Having shown that the members of $\CD {G}$ are subnormal, one continues by asking about the Chermak-Delgado lattice of a direct product. Before proceeding, though, we introduce a subset of the Chermak-Delgado lattice.  

In \cite{glauberman2006} Glauberman defines the notion of a {\it centrally large subgroup} and shows, among other things, that a subgroup $U$ is centrally large exactly when $U \in \CD G$ and $\Z U = C_G (U)$.  We denote the set of centrally large subgroups of $G$ by $\CL G$.  Note that $\CL G$ is closed under joins and contains the largest element in $\CD {G}$.  

In addition to describing $\CD {G \times H}$ for finite groups $G$ and $H$, we also describe $\CL {G \times H}$.  We utilize the following basic fact about centralizers in direct products, the proof of which follows directly from the mechanics of conjugation in a direct product.

\begin{lem}\label{directcentralizers} Let $G$ and $H$ be groups.  If $U \leq G \times H$ then
$C_{G \times H} (U) = C_G (\pi_G(U)) \times C_H (\pi_H(U))$.
\end{lem}

\begin{thm} \label{direct} For any finite groups $G$ and $H$, the lattices $\CD {G \times H}$ and $\CD {G} \times \CD {H}$ are equal and $\CL {G \times H} = \CL {G} \times \CL {H}$.
\end{thm}

\begin{proof}
Let $U \leq G \times H$.  We have the following inequality, with the second step due to Lemma~\ref{directcentralizers}.
$$
\begin{array}{rl}
m_{G \times H} (U) &= \ord U \ord {C_{G \times H} (U)}\\
&= \ord U \ord {C_G (\pi_G(U)) \times C_G(\pi_H(U))}\\
&\leq \ord {\pi_G(U) \times \pi_H(U)} \ord {C_G (\pi_G(U)) \times C_G(\pi_H(U))}\\
&\leq \ord {\pi_G(U)} \ord {C_G(\pi_G(U))} \ord {\pi_H(U)} \ord {C_H(\pi_H(U))}\\
&\leq m_G(\pi_G(U)) m_H(\pi_H(U))\\
\end{array}
$$
Equality occurs {\it exactly} when $U = \pi_G(U) \times \pi_H(U)$. Therefore, the subgroups of $G \times H$ with maximal measure are exactly those direct products $X \times Y$ where $X \in \CD {G}$ and $Y \in \CD {H}$. This gives $\CD {G \times H} = \CD {G} \times \CD {H}$. 

Now suppose that $U \in \CL {G \times H}$; then $U \in \CD {G \times H}$. For $X \in \{G, H \}$ we have $\pi_X (U) \in \CD X$.  Let $g \in C_G(\pi_G(U))$.  The element $(g,1)$ centralizes $U$, hence its projection $g$ is in $C_G(\pi_G(U))$.  Therefore $C_G(\pi_G(U)) \leq \pi_G (U)$; we can similarly prove the same with respect to $H$.  Hence $\pi_X (U) \in \CL {X}$ for $X = G, H$.

Alternatively assume that $\pi_X (U) \in \CL {X}$ for $X \in \{G, H\}$.  Then $\pi_X (U) \in \CD {X}$ and hence $U \in \CD {G \times H}$.  Moreover, $\Z {U \cap X} = C_X (\pi_X (U))$ for both values of $X$, thus $\Z U = C_{G \times H} (U)$ after using Lemma~\ref{directcentralizers}.  Therefore $U \in \CL {G \times H}$. Hence $\CL {G \times H} = \CL {G} \times \CL {H}$, as desired.\end{proof}

\section{Wreath Products}

This section discusses our attempts to describe the Chermak-Delgado lattice of a wreath product. As a byproduct of our efforts we show that every finite 2-group $G$ can be embedded in a finite 2-group $E$ such that $E \in \CD E$, while noting that there are 2-groups $E$ such that $E \not \in \CD E$.

Let $G$ and $H$ be finite groups with $H$ a permutation group of degree $n$ acting on a set $\Omega$.  The {\it wreath product of $G$ by $H$}, denoted $G \wr H$, is the semidirect product $B \rtimes H$ where $B = G^\Omega$ is the group of all functions $f: \Omega \rightarrow G$ under point-wise multiplication.  The subgroup $B$ is referred to as the {\it base} of $W$.  If $h \in H$ and $f \in B$ then
$$
f^h (\omega) = f(\omega h^{-1})
$$
for $\omega \in \Omega$.  

We focus on wreath products $G \wr H$ where $H \iso C_n$ for some positive integer $n$, so $\Omega = \{ 1, 2, \dots, n \}$.  As in the case with direct products, we start by examining centralizers.

\begin{prop} \label{centralizers} Let $G$ be a non-trivial group and set $W = G \wr C_n$ where $C_n = \langle \sigma \rangle$ is the cyclic group of order $n$.  Let $B$ be the base group of $W$.  If, for some $f \in B$, the element $f\sigma \in W - B$ commutes with an element $b \in B$ then 
$$
b(i) = b(1)^{f(1)f(2) \cdots f(i-1)} \quad \textrm{for } 1 < i \leq n.
$$
Thus all $b(i)$ are in some orbit of $\langle f(1), f(2), \dots, f(n) \rangle$.  Furthermore 
$$b(1) \in C_G (f(1) f(2) \cdots f(n))$$
and hence $\pi_1 (C_B(f \sigma)) \iso C_G (f(1)f(2) \cdots f(n))$. \end{prop}

\begin{proof} Suppose that for $f \in B$ and $f \sigma \in W - B$ commutes with some $b \in B$.  Notice: 
$$
f \sigma b = b f \sigma \iff b^{\sigma^{-1}} = b^f.
$$
In particular, for $1 < i \leq n$ we have
$$
b^{\sigma^{-1}} (i) = b(i \sigma) = b^f (i) = b(i)^{f (i)}.
$$
Plugging in a few values for $i$, we see $b(2) = b(1)^{f(1)}$ and
$$
b(3) = b(2)^{f(2)} = (b(1)^{f(1)})^{f(2)} = b(1)^{f(1)f(2)}.
$$
Continuing in this way, one concludes that $b(i) = b(1)^{f(1)f(2) \cdots f(i-1)}$ for all $i$ with $1 < i \leq n$.  

In fact, since $n \sigma = 1$ we also see that
$$
b(1) = b(n \sigma) = b(n)^{f(n)} = b(1)^{f(1)f(2) \cdots f(n)},
$$
hence $b(1)$ commutes with $f(1)f(2) \cdots f(n)$ in $G$ and 
$$\pi_1 (C_B(b\sigma)) \iso C_G (f(1)f(2) \cdots f(n)),$$
as described in the statement of the proposition.
\end{proof}

This proposition can generalize straightforwardly to more general $H$ but since the notation quickly becomes cumbersome and we do not apply such a generalization here, we refer the reader to \cite{diss}. Proposition~\ref{centralizers} is enough to establish some facts about $C_W (B)$ and $\Z W$, allowing us to better calculate $m_W (B)$ and $m_W (W)$.  

\begin{prop} \label{basecent} Let $W = G \wr C_n$ with $G$ a non-trivial group and base group $B$.  The centralizer in $W$ of $B$ is $\Z B$; consequently $m_W(B) = m_B(B) = {\ord G}^n {\ord {\Z G}}^n$. \end{prop}

\begin{proof} Set $C_n = \langle \sigma \rangle$.  Suppose an element $z \in W$ centralizes $B$.  If there exists $f \in B$ such that $z = f \sigma \not \in B$ then Proposition~\ref{centralizers} redefines structure of $B$, namely telling us that $B \iso C_G (f(1) f(2) \cdots f(n))$.  Yet by its definition $B$ cannot be isomorphic to a subgroup of $G$.  Hence $z$ must be an element of $B$, yielding $C_W(B) \leq B$.  Thus $C_W (B) = \Z B$. Therefore $m_W(B) = \ord {B} \ord {\Z B} = {\ord G}^n {\ord {\Z G}}^n$ as claimed. \end{proof}

Combining Proposition~\ref{basecent} with \cite[Exercise 3A.9]{Isaacs}, which states that elements commuting with the generator of $C_n$ must be diagonal, we have the following description of $\Z W$.

\begin{prop} \label{center} Let $G$ be a non-trivial group, set $W = G \wr C_n$, and let $B$ represent the base of $W$.  The center of $W$ is equal to the diagonal of $\Z B$ and consequently $m_W (W) = n {\ord G}^n \ord {\Z G}$. \end{prop}

The next proposition is a straightforward consequence of our calculations in Propositions~\ref{basecent} and \ref{center}.  The result implies that even when $G \in \CD G$ and $H \in \CD H$, the wreath product $W = G \wr H$ need not be a member of $\CD W$.  

\begin{prop} \label{p>2} Let $G$ be a group and let $W = G \wr C_n$ for an integer $n \geq 2$.  If $\ord {\Z G} \geq 2$ or $n > 2$ then $W \not \in \CD W$. \end{prop}

\begin{proof}  Let $z = \ord {\Z G}$.  We first calculate the measures of $W$ and $B$ using Propositions~\ref{center} and \ref{basecent}:
$$
m_W(W) = n {\ord G}^n \ord {\Z G} = n {\ord G}^n \cdot z \textrm{ and}
$$
$$
m_W(B) = {\ord G}^n {\ord {\Z G}}^n = {\ord G}^n \cdot z^n.
$$
Thus $m_W(W) < m_W(B)$ if and only if $z > n^{1/{n-1}}$.  One easily confirms that this latter expression is strictly decreasing for integers $n \geq 2$.  When $n = z = 2$ or when $\ord {\Z G} = 1$ we have $m_W(W) \geq m_W(B)$.  Otherwise, though, $m_W(W) < m_W(B)$ and thus $W \not \in \CD W$. \end{proof}

Observe from the proof of Proposition~\ref{p>2} that when $\ord {\Z G} = n = 2$ then $m_W (W) = m_W (B)$. We saw an example of this situation, $D_8$, where $W \in \CD W$. Therefore, in light of $\CD {D_8}$ and Proposition~\ref{p>2}, we are interested two questions:  

\begin{enumerate}
\item If $\ord {\Z G} = 2$ will $W = G \wr C_2$ be a member of $\CD{W}$?  
\item If $W = G \wr C_n$ with $\ord {\Z G} >2$ or $n > 2$ will $\CD {W} = \CD {B}$? 
\end{enumerate}

\noindent
In the remainder of the section we address both of these questions.  Let us first note that if $G$ is not in its own Chermak-Delgado lattice then $W$ need not be in $\CD {W}$.  The first nonabelian group $G$ with $\Z G \iso C_2$ and $G \not \in \CD {G}$ is $D_{12}$, the dihedral group of order 12.  

\begin{example}  Let $G = D_{12}$. First we show that $G \not \in \CD {G}$. Let $r$ be an element of order 6; then $\langle r \rangle = C_G (\langle r \rangle)$.  Hence $m_{G} (\langle r \rangle) = 6^2 = 36$.  Yet $m_{G} (G) = 12 \cdot 2 = 24$, so $G \not \in \CD {G}$.  

Let $W = G \wr C_2$; then $m_W (W) = 6^2 \cdot 2^4$.  Let $U$ be the subgroup of the base of $W$ isomorphic to $\langle r \rangle \times \langle r \rangle$.  Observe that $U \leq C_W (U)$ and therefore $m_W (U) \geq {\ord U}^2 = 6^4$.  Since $m_W (U) > m_W (W)$, we know that $W \not \in \CD {W}$. \end{example}

Thus, with regards to Question 1, we show that if $\ord {\Z G} = 2$ {\it and} $G \in \CD {G}$ then $W = G \wr C_2$ is in $\CD {W}$ and $\CD {B} \leq \CD {W}$ as lattices.  To attain this answer and to address Question 2, we examine $\CD {W}$ by considering $m_W (U)$ for $U \in \CD {W}$.  There are four cases, depending upon whether or not $U \leq B$ or $C_W (U) \leq B$.  The next lemma describes a reduction in calculuating the order of $U$; it's a direct consequence of the Isomorphism Theorems \cite[Theorem 3.18]{DF}.

\begin{lem}\label{index} Let $G$ be a non-trivial group, $W = G \wr C_p$ for some prime $p$, and $B$ be the base of $W$.  If $U \leq W$ then $\ind {U}{B \cap U} = \begin{cases}
1 & \textrm{if } U \leq B\\
p & \textrm{if } U \not \leq B\\
\end{cases}$.\end{lem}

We use Lemma~\ref{index} in the proof of the following result, the key observation for calculating the Chermak-Delgado measure of a subgroup in a wreath product.

\begin{prop}\label{orders} Let $G$ be a non-trivial group and let $W = G \wr C_p$ for a prime $p$.  Suppose $B$ is the base of $W$ and let $U \leq W$.  

\begin{enumerate}
\item If $U \leq B$ and $C_W (U) \not \leq B$ then $\ord U = \ord {\pi_1 (U)}$ and $\ord {C_W(U)} = p {\ord {C_G(\pi_1(U))}}^p$.
\item If $U \not \leq B$ and $C_W(U) \not \leq B$ then $\ord U = p \ord {\pi_1 (U \cap B)}$ and $\ord {C_W(U)} = p \ord {\pi_1 (C_B(U))}$.
\end{enumerate} \end{prop}

\begin{proof}Let $U \leq B$ and suppose $C_W (U) \not \leq B$. After applying Lemma~\ref{index} to $C_W (U)$, we see that $\ord {C_W (U)} = p \ord {C_B (U)}$; moreover, $C_W (U) / C_B (U) \iso W / B \iso \langle \sigma \rangle$ and there must exist $f \in B$ such that $C_W (U) = C_B (U) \langle f \sigma \rangle$.  

Proposition~\ref{centralizers} then applies to $U \leq B$ and $f \sigma \in C_W (U)$, so that if $u \in U$ there exists a $g \in C_G (f(1) f(2) \cdots f(p-1))$ with
$$
u(i) = g^{f(1)f(2) \cdots f(i-1)} \textrm{ for each } i \in \Omega.
$$
Thus $\pi_i(U) = (\pi_1(U))^{f(1)f(2) \cdots f(i-1)}$ for $2 \leq i \leq p$.  Therefore $U$ is a ``diagonal-type'' subgroup and $\ord U = \ord {\pi_1(U)}$, as claimed.  Moreover, the description of $U$ from Proposition~\ref{centralizers} implies that $\ord U = \ord {\pi_1(U)} \leq \ord {C_G (f(1) f(2) \cdots f(p-1))}$.

Lemma~\ref{directcentralizers} states that $C_B (U) = \prod\limits_{i=1}^p C_G (\pi_i (U))$.  Given the structure of $U$, we can establish $\pi_1(U) \iso \pi_2(U)^{f(1)}$ and, similarly, $\pi_i (U) = (\pi_1(U))^{f(1)f(2) \cdots f(i-1)}$ for all $i$ with $3 \leq i \leq p$.  Therefore $C_G (\pi_1(U)) \iso C_G (\pi_i(U))$ for all $i$ with $2 \leq i \leq p$, and $\ord {C_B (U)} = \ord {C_G (\pi_1(U))}^p$.

Now suppose that neither $U$ nor $C_W(U)$ is a subgroup of $B$.  Then Lemma~\ref{index} tells us that $\ord U = p \ord {U \cap B}$ and $\ord {C_W(U)} = p \ord {C_B(U)}$.  Yet $U \cap B \leq B$ and $C_W(U \cap B)$ contains $C_W(U)$, and hence $C_W (U \cap B) \not \leq B$.  Applying part (1) to $U \cap B$ we have $\ord {U \cap B} = \ord {\pi_1 (U \cap B)}$.  Thus 
$$
\ord U = p \ord {\pi_1 (U \cap B)},
$$
as desired.  

Let $X = C_W (U)$.  Note that $X \not \leq B$ and we established $\ord X = p \ord {X \cap B}$.  Since $U \leq C_W(X)$ we know $C_W (X) \not \leq B$.  Apply the argument of the last paragraph to $X$; thus $\ord X = p \ord {\pi_1 (X \cap B)}$.  Since $X \cap B = C_B(U)$, we therefore have shown that $\ord {C_W(U)} = p \ord {\pi_1 (C_B(U))}$. \end{proof}

\begin{thm} \label{wreath} Let $G \in \CD {G}$ and suppose $\ord {\Z G} = 2$. Let $W = G \wr C_2$.  The group $W$ is a member of $\CD {W}$ and $\CD {B} \leq \CD {W}$, as lattices. \end{thm}

\begin{proof} First we calculate the measures of $W$ and $B$, using Propositions~\ref{center} and \ref{basecent}.  This gives:
$$
m_W(W) = 2^2 {\ord G}^2 = m_W (B) = m_B (B).
$$
We will show for all $U \in \CD {W}$ that $m_W(U) \leq 2^2 {\ord G}^2$, thus determining that $W$ has maximal measure and implying $W \in \CD {W}$.  To do this we will first consider $U \leq B$ and then turn our attention to $U \not \leq B$.  

If $C_W (U) \leq B$ then $C_W (U) = C_B(U)$.  Thus $m_W(U) = m_B(U)$.  Since $G \in \CD {G}$, we know by Theorem~\ref{direct} that $B \in \CD {B}$.  Therefore 
$$
m_W (U) = m_B(U) \leq m_B(B) = m_W(W).
$$
If, on the other hand, when $U \leq B$ we also have $C_W (U) \not \leq B$ then Proposition~\ref{orders} yields that $\ord U = \ord {\pi_1(U)}$.  This, together with the information about the centralizer of $U$ from Proposition~\ref{orders}, yields
$$
\begin{array}{rl}
m_W(U) &= \ord U \ord {C_W(U)}\\
& \leq \ord {\pi_1 (U)} \cdot 2 \cdot \ord {C_G (\pi_1(U))}^2\\
& \leq 2 \cdot m_G (\pi_1 (U)) \cdot \ord {C_G (\pi_1(U))}.\\
\end{array}
$$
Yet $G \in \CD {G}$; thus $m_G (\pi_1 (U))$ is less than $2 \ord G$.  Also, the centralizer of $\pi_1 (U)$ clearly has order no more than $\ord G$.  This allows $m_W(U) \leq 2^2 {\ord G}^2$. Therefore if $U \leq B$ then $m_W(U) \leq m_W(W)$. 

Now suppose that $U \not \leq B$. If $C_W(U) \leq B$ then we know already that
$$
m_W(C_W (U)) \leq m_W(W),
$$
by the preceding paragraphs.  Yet $U \in \CD {W}$, so $m_W(U) = m_W (C_W(U))$.  Hence we need only examine the case where $C_W (U) \not \leq B$.  

In this case, Proposition~\ref{orders} tells us that $\ord U = 2 \ord {\pi_1 (U \cap B)}$ and $\ord {C_W(U)} \leq 2 \ord {\pi_1(C_B (U))}$.  Notice that $C_B (U) \leq C_B (U \cap B)$, and therefore $\pi_1 (C_B (U)) \leq \pi_1 (C_B (U \cap B))$.  It's a straightforward argument to show that $\pi_1 (C_B (U \cap B)) \leq C_G (\pi_1 (U \cap B))$.  Therefore
$$
\begin{array}{rl}
m_W(U) &= \ord U \ord {C_W(U)}\\
&\leq 2 \cdot \ord {\pi_1(U \cap B)} \cdot 2 \cdot \ord {\pi_1 (C_G (U \cap B))}\\
&\leq 2^2 \cdot \ord {\pi_1(U \cap B)} \ord {C_G (\pi_1 (U \cap B))}\\
&\leq 2^2 m_G (\pi_1 (U \cap B)).\\
\end{array}
$$
Yet $G \in \CD {G}$, so we can conclude that
$$
m_W (U) \leq 2^3 \ord G.
$$
Since $\ord G \geq 2$, this yields the desired result for $U \not \leq B$.

To finish the proof, let $U \in \CD {B}$.  Then $m_B (U) = m_B (B)$, yet we've established that this latter quantity equals $\M {W}$.  Hence $U \in \CD {W}$, as well.\end{proof}

The proof of Theorem~\ref{wreath} establishes that when $W$ is as described then $\CD {W}$ contains at least $\CD {B}$ and new maximal and minimal elements ($W$ and $\Z W$, respectively). There may even be other elements of $\CD {W}$ that are not in $\CD {B}$.  And, as a corollary of Theorem~\ref{wreath}, we have the result mentioned at the start of the section.

\begin{cor} \label{embed} If $G$ is a 2-group then there exists a 2-group $E$ with $E \in \CD {E}$ such that $G$ can be embedded as a subgroup of $E$. \end{cor}

\begin{proof}  The group $G$ can be embedded as a subgroup of $S_n$ for some $n$.  Let $E$ be the Sylow 2-subgroup of $S_n$ that contains $G$.  Recall that $E$ is a direct product whose factors are iterated wreath products of $C_2$. Each of the iterated wreath products is contained in its Chermak-Delgado lattice, by Theorem~\ref{wreath}.  Thus $E \in \CD {E}$ by Theorem~\ref{direct}. \end{proof}  

Corollary~\ref{embed} is not trivial, in the sense that there are 2-groups which are not in their own Chermak-Delgado lattice.  In fact, there are 2-groups $G$ with $\Z G = 2$ such that $G \notin \CD {G}$.  We provide one example here.

\begin{example} Let $G$ be the Sylow 2-subgroup of the general linear group of $n \times n$ matrices with entries in the field of order 2.  It is known that $G$ is isomorphic to the group of upper triangular matrices over the field of order 2 and that $\ord {\Z G} = 2$.  

Let $A$ be an abelian subgroup of maximal rank in $G$.  In \cite{thwaites} it's shown that $\ord A = 2^{xy}$ where $x$ is the greatest integer less than or equal to $\frac{n}{2}$ and $y$ is the smallest integer greater than or equal to $\frac{n}{2}$; thus $m_G (A) \geq 2^{2xy}$.  On the other hand, $m_G (G) = 2^{n(n-1)/2} \cdot 2$.  When $n = 5$, we have $x = 2$ and $y = 3$ so it's easy to see that $m_G(A) > m_G (G)$.\end{example}

Recall Question 2:  If $W = G \wr C_n$ where $\ord {\Z G} > 2$ or $n > 2$, will $\CD {W} = \CD {B}$?  We address this question only in the case where $n$ is a prime number.  The techniques to give an affirmative answer for this restricted question are along the same lines as what we have done so far in this section.  We begin with a lemma.

\begin{lem}\label{mmm}  Let $p$ be a prime number, $G$ be a group with $\Z G > 1$, and $W = G \wr C_p$ with base $B$.  If $W \not \iso D_8$ then for every $U \in \CD W$ either $U \leq B$ or $C_W (U) \leq B$. \end{lem}

\begin{proof}  We prove the contrapositive of the lemma, supposing that there exists $U \in \CD {W}$ with $U \not \leq B$ and $C_W (U) \not \leq B$.  Then Proposition~\ref{orders} yields:
$$
\begin{array}{rl}
\M {W} 	& = \ord U \ord {C_W (U)}\\
		& \leq p^2 \ord {\pi_1 (U \cap B)} \ord {\pi_1 (C_B (U))}\\
		& \leq p^2 \ord {\pi_1 (U \cap B)} \ord {\pi_1 (C_B (U \cap B))}\\
		& \leq \M {G} \cdot p^2.\\
\end{array}
$$
At the same time, though, we know that $\M {B} = (\M {G})^p \leq \M {W}$.  Thus
$$
{(\M G)}^p \leq \M G \cdot p^2.
$$
The usual algebra tactics allow us to rearrange the inequality:  $\M G \leq (p^2)^{1/(p-1)}$.  Yet this last expression is the square of a function that is strictly decreasing on integers $n \geq 2$; hence its maximum value is when $p = 2$.  Additionally, $m_G (G) \leq \M G$ and therefore 
$$
\ord G \ord {\Z G} \leq (p^2)^{1/(p-1)} \leq 4.
$$
Since $\ord {\Z G} \geq 2$, the above can only occur when $\ord {\Z G} = \ord G = p = 2$.  In this case, though, $W \iso D_8$. \end{proof}

\begin{thm}\label{nnn} Let $p$ be a prime, $G$ be a group with $\Z G > 1$, and $W = G \wr C_p$ with base $B$.  If $\ord {\Z G} > 2$ or $p > 2$ then for every $U \in \CD {W}$ both $U \leq B$ and $C_W(U) \leq B$. Thus in this case $\CD W = \CD B$ and similarly, $\CL W = \CL B$. \end{thm}

\begin{proof}  We prove the contrapositive of the theorem.  Assume that there exists $U \in \CD W$ and at least one of $U$ or $C_W(U)$ is not a subgroup of $B$.  If both $U$ and $C_W (U)$ are not subgroups of $B$ then Lemma~\ref{mmm} tells us that $W \iso D_8$.  In this case $\ord {\Z G} = p = 2$, so the theorem holds.  

Suppose exactly {\it one} of $U$ or $C_W (U)$ is not a subgroup of $B$.  As $U \in \CD W$, we know $U = C_W (C_W (U))$, therefore we may assume, without loss of generality, that $U \leq B$ and $C_W (U) \not \leq B$.

Proposition~\ref{orders} implies
$$
\begin{array}{rl}
\M {W} 	& = \ord U \ord {C_W(U)}\\
   		& = \ord {\pi_1 (U)} \ord {C_G(\pi_1(U))}^p \cdot p\\
		& = \M {G} \cdot \ord {C_G (\pi_1(U))}^{p-1} \cdot p.\\
\end{array}
$$
Again we know that $\M B = (\M G)^p \leq \M W$ and additionally $\ord {C_G (\pi_1(U))} \leq G$.  Hence
$$
(\M G)^p \leq \M G \cdot {\ord G}^{p-1} \cdot p \textrm{ and therefore } \M G \leq \ord G \cdot p^{1/(p-1)}.
$$
Then $\ord G \ord {\Z G} \leq \M {G} \leq \ord G \cdot p^{1/(p-1)}$ and therefore $\ord {\Z G} \leq p^{1/(p-1)}$.  This familiar expression is strictly decreasing on integers $n \geq 2$, as before.  Therefore $\ord {\Z G} \leq 2$ and, given the hypotheses of the theorem, $\ord {\Z G} = p = 2$.  

Therefore if $\ord {\Z G} > 2$ or $p > 2$ then for every $U \in \CD {W}$ we know $U \leq B$ and $C_W(U) = C_B(U)$.  Thus $m_W(U) = m_B(U)$.  It is always true that $\M {B} \leq \M {W}$ so in this case $U \in \CD {B}$ and $\CD {W} \leq \CD {B}$.  That then implies that $\M {W} = \M {B}$, since $C_W(U) \leq C_B(U)$ for any $U \leq B$.  Hence if $U \in \CD {B}$ then $U \in \CD {W}$, too.  Thus $\CD {W} = \CD {U}$.

Let $U \in \CL {W}$.  Then $U \in \CD {W}$ and $C_W (U) = \Z U$.  By the arguments earlier, $U \in \CD {B}$.  Additionally $C_B (U) \leq C_W (U)$, so we conclude $U \in CL {B}$, giving $\CL {W} \subseteq \CL {B}$.  Suppose that $U \in \CL {B}$; hence $U \in \CD {B}$ and $C_B (U) = \Z U$.  Then we can conclude $U \in \CD {W}$, but then the preceeding paragraph gives $C_W (U) \leq C_B (U)$.  Therefore $U \in \CL {W}$ and $\CL {B} = \CL {W}$. \end{proof}

\section*{Acknowledgements}

We thank both Peter Hauck, University of T\"ubingen, and Marcin Mazur, Binghamton University, for their suggestions regarding the proof of Theorem~\ref{subnormal}. We also thank Joseph Bohannon for his aid in computing the example following Theorem~\ref{subnormal}.

\bibliographystyle{srtnumbered}
\bibliography{references}

\end{document}